\documentclass{article}
\usepackage{latexsym}
\usepackage[varumlaut]{yfonts}
\usepackage{epsfig}
\usepackage{amsmath}
\usepackage{amssymb}
\usepackage{amsthm}
\usepackage{bm}
\newtheorem{theorem}{Theorem}[section]

\newtheorem{lemma}[theorem]{Lemma}

\theoremstyle{remark}

\theoremstyle{definition}

\DeclareMathOperator\Aut{Aut}

\DeclareMathOperator\diag{diag}

\DeclareMathOperator\sign{sgn}

\DeclareMathOperator\spa{span}

\begin{document}

\title{Analytic formulas for complete hyperbolic affine spheres}

\author{Roland Hildebrand \thanks{%
LJK, Universit\'e Grenoble 1 / CNRS, 51 rue des Math\'ematiques, BP53, 38041 Grenoble cedex 09, France
({\tt roland.hildebrand@imag.fr}, tel.\ +33(0)476-63-5714, fax +33(0)476-63-1263).}}

\maketitle

\begin{abstract}
We classify all regular three-dimensional convex cones which possess an automorphism group of dimension at least two, and provide analytic expressions for the complete hyperbolic affine spheres which are asymptotic to the boundaries of these cones. The affine spheres are represented by explicit hypersurface immersions into three-dimensional real space. The generic member of the family of immersions is given by elliptic integrals.
\end{abstract}



\section{Introduction}

When studying a class of mathematical objects, it is always advantageous to have examples at hand which possess an explicit analytic description. Affine spheres, a class of hypersurfaces in affine real spaces studied in affine differential geometry, are described by solutions of a Monge-Amp\`ere equation, a fully nonlinear elliptic partial differential equation (PDE). Analytic solutions of nonlinear PDEs being generally scarce, only a few explicit descriptions of affine spheres are known. On the other hand, it is well-known that the class of affine spheres is large. The Calabi conjecture \cite[p.22]{Calabi72}, which was proven around the 1980s by the efforts of many authors (see \cite[Section 2]{LSZ} for a synthesis of the proof), states that every regular convex cone possesses a unique homothetic foliation by complete hyperbolic affine spheres which are asymptotic to the boundary of the cone, and all complete hyperbolic affine spheres can be obtained in this way. Here a convex cone is called {\it regular} if it is closed, has non-empty interior, and does not contain a line.

However, all complete hyperbolic affine spheres for which an explicit analytic description is known are asymptotic to the boundary of homogeneous cones. As such, they are level sets of the characteristic function of the cone, and can actually be obtained by easier means than solving the PDE. No analytic example of an affine sphere is known which is not homogeneous. This situation cannot be considered to be satisfactorily, especially when one wants to study properties of affine spheres that become non-trivial only in a non-homogeneous setting. Examples of such properties are variations of different norms of the curvature or the cubic form over the affine sphere.

In this paper we provide an analytic description of families of non-homogeneous complete hyperbolic affine spheres. The key idea is that if a cone $K$ possesses an automorphism group with orbits of codimension 1, then the Monge-Amp\`ere equation can be reduced to an ordinary differential equation (ODE). This ODE is then much easier to solve that the original PDE, and one may hope to obtain an analytic solution. The simplest case is that of a regular convex three-dimensional cone with an automorphism group of dimension at least 2. In this contribution we classify all such cones, which we may call {\it semi-homogeneous}. We deduce and solve the ODEs describing the hyperbolic affine spheres which are asymptotic to the boundary of these cones. In general, the solution will be given by elliptic functions.

The remainder of the paper is structured as follows. In the next section we introduce the notion of an affine sphere and describe a convenient representation of affine spheres as level sets of a logarithmically homogeneous function. In Section \ref{sec:semi_hom} we define the semi-homogeneous cones and provide their classification in three dimensions. In Section \ref{sec:affine_spheres} we compute the complete hyperbolic affine spheres which are asymptotic to the boundary of the three-dimensional semi-homogeneous cones.

\section{Affine spheres and the canonical potential}

Equiaffine differential geometry studies the properties of submanifolds in affine space $\mathbb A^n$ which are invariant under the group of volume preserving affine transformations. On a non-degenerate hypersurface $M \subset \mathbb A^n$ there exists a unique, up to a sign change, affinely invariant transversal vector field, the {\it affine normal field} \cite[p.42]{NomizuSasaki}. The affine lines defined by this vector field are called {\it affine normals}. If we consider the affine normal field as a unit normal vector field, then the second fundamental form of the hypersurface defines a quadratic form on the tangent bundle of the hypersurface, the {\it affine fundamental form}. If the hypersurface is convex, then the affine fundamental form is definite and defines a Riemannian metric on the hypersurface, the {\it affine metric}.

If all affine normals meet in one point, then the hypersurface is called a {\it proper affine sphere} \cite[p.43]{NomizuSasaki}. The intersection point of the affine normals is called the {\it centre} of the affine sphere. If the proper affine sphere is convex and separates its convex hull from its centre, then it is called {\it hyperbolic}. Hyperbolic affine spheres which are complete in the affine metric are the subject of the Calabi conjecture \cite[p.22]{Calabi72}. This conjecture states that every complete hyperbolic affine sphere is asymptotic to the boundary of a regular convex cone $K \subset \mathbb A^n$, with the vertex of the cone being located at the centre of the affine sphere, and that every regular convex cone gives rise to a unique 1-parametric family of complete hyperbolic affine spheres, related by homothety with respect to the vertex of the cone, which are asymptotic to the boundary of the cone. This provides a full classification of complete hyperbolic affine spheres.

In the sequel it will be convenient to consider the cones $K$ as subsets of the real vector space $\mathbb R^n$, with the vertex of $K$ located at the origin. The complete hyperbolic affine spheres which are asymptotic to the boundary of $K$ will be described as the level sets of the solution $F: K^o \to \mathbb R$ of the following Monge-Amp\`ere equation \cite[eq.(4.1), p.359]{ChengYau82}
\begin{equation} \label{MAequation}
\det F'' = e^{2F},\qquad F|_{\partial K} = +\infty,\qquad F'' \succ 0.
\end{equation}
By \cite{ChengYau82} the solution of \eqref{MAequation} exists for every regular convex cone $K$. From \cite[Prop.\ 5.5, p.528]{ChengYau80} it follows that this solution is unique, see \cite{Hildebrand12b} for more details. A more comprehensive review of the literature and an independent derivation can be found in \cite{Fox12a}, where the solution of \eqref{MAequation} was termed {\it canonical potential} of the cone $K$. Note that the determinant in the left-hand side of the equation implies the availability of a volume form on $\mathbb R^n$. The solution will depend on the choice of this volume form. Namely, multiplication of the volume form by $\alpha > 0$ leads to the transformation $F \mapsto F + \frac12\log\alpha$. This transformation does not change the level sets of $F$, however, and hence the choice of the volume form will be irrelevant for our purposes.

By making the substitution $\tilde x = \alpha x$, $\alpha > 0$, and comparing the solution of \eqref{MAequation} in both coordinate systems $x$ and $\tilde x$ on $\mathbb R^n$, one easily sees that the function $F$ is {\it logarithmically homogeneous} of degree $-n$, i.e., for all $\alpha > 0$ and all $x \in K^o$ we have
\begin{equation} \label{log_homogeneous}
F(\alpha x) = -n\log\alpha + F(x),
\end{equation}
where $n$ is the dimension of $K$. That the solution of \eqref{MAequation} on regular convex cones is logarithmically homogeneous was already mentioned in \cite[pp.426--427]{Loftin02a} without a detailed proof. We have the following result.

{\theorem (\cite[Theorem 3]{Loftin02a}, see also \cite[Appendix A]{Sasaki85}) Let $K \subset \mathbb R^n$ be a regular convex cone. Then the complete hyperbolic affine spheres which are asymptotic to the boundary of $K$ are exactly the level sets of the solution $F: K^o \to \mathbb R$ of equation \eqref{MAequation}. }

Let $\Aut K$ denote the automorphism group of the cone $K$. Even for simple cones $K$ such as polyhedral cones except the orthant, an analytic expression for the canonical potential is not known. However, the equiaffine invariance of equation \eqref{MAequation} implies that the canonical potential $F$ is invariant under unimodular automorphisms $g \in \Aut K$ of the cone $K$, see \cite{Hildebrand12b} or \cite{Fox12a}. Thus when solving \eqref{MAequation}, the dimension of the problem can be effectively reduced by the generic dimension of the orbits of $\Aut K$. If there are orbits of dimension $n-1$, then the PDE \eqref{MAequation} reduces to an ODE. In the next two sections we will apply this approach to three-dimensional regular convex cones possessing orbits of dimension at least two.

\section{Semi-homogeneous cones in $\mathbb R^3$} \label{sec:semi_hom}

In this section we classify the three-dimensional regular convex cones with automorphism group of dimension at least 2. Let $K \in \mathbb R^3$ be a regular convex cone. Then $\dim\Aut K \geq 2$ if and only if there exists an element $A$ in the Lie algebra $\mathfrak{aut}K$ of $\Aut K$ which is not proportional to the identity matrix. Our first step will be to derive a canonical form for $A$, which is given by the following lemma.

{\lemma \label{can_form} Let $K \in \mathbb R^3$ be a regular convex cone such that $\dim\Aut K \geq 2$. Then there exists a regular convex cone $\tilde K \in \mathbb R^3$ and an element $A$ in the Lie algebra $\mathfrak{aut}\tilde K$ of $\Aut \tilde K$ such that $\tilde K$ is isomorphic to $K$ and $A$ is given by one of the matrices
\[ \begin{pmatrix} 0 & 1 & 0 \\ 0 & 0 & 1 \\ 0 & 0 & 0 \end{pmatrix},\ \begin{pmatrix} 0 & 1 & 0 \\ 0 & 0 & 0 \\ 0 & 0 & 0 \end{pmatrix},\ \begin{pmatrix} 1 & 1 & 0 \\ 0 & 1 & 0 \\ 0 & 0 & 0 \end{pmatrix},\ \begin{pmatrix} 1 & 0 & 0 \\ 0 & -\mu & 0 \\ 0 & 0 & \mu-1 \end{pmatrix},\ \begin{pmatrix} \mu & 1 & 0 \\ -1 & \mu & 0 \\ 0 & 0 & 0 \end{pmatrix},
\]
where in the fourth case $\mu \in [0,\frac12]$, and in the last case $\mu \geq 0$. }

\begin{proof}
Let $S \in GL(3,\mathbb R)$ be an arbitrary automorphism of $\mathbb R^3$, and let $\tilde K = S[K]$ be the image of $K$. Then $S$ induces an isomorphism between the Lie algebras $\mathfrak{aut}K,\mathfrak{aut}\tilde K$ given by the similarity transformation $A \mapsto S^{-1}AS$. We may hence assume without restriction of generality that there exists $A \in \mathfrak{aut}\tilde K$ which is not proportional to the identity matrix and is in real Jordan form.

Since the Lie algebra $\mathfrak{aut}\tilde K$ is a linear space, we may multiply $A$ with an arbitrary nonzero constant. Moreover, since homotheties are always automorphisms of cones, and hence the identity matrix $I$ is always an element of $\mathfrak{aut}\tilde K$, we may add an arbitrary multiple of the identity matrix to $A$. Thus we may without restriction of generality replace $A$ with the linear combination $\alpha A + \beta I$, where $\alpha \not= 0$.

If $A$ has three equal real eigenvalues, then it might be brought to one of the first two canonical forms in the lemma, depending on whether the number of Jordan cells is 1 or 2. If there are two distinct real eigenvalues, but the matrix is not diagonalizable, then it might be brought to the third canonical form. This can be done by first normalizing the single eigenvalue to zero by adding a multiple of the identity matrix, then normalizing the double eigenvalue to 1 by multiplying with an appropriate constant, and finally passing to the real Jordan form. If all eigenvalues of $A$ are real and the matrix is diagonalizable, then it might be brought to the fourth canonical form. To see this, first normalize the trace to zero by adding an appropriate multiple of the identity matrix, then normalize the eigenvalue with the largest absolute value to 1 by multiplying with an appropriate constant, and finally arrange the eigenvalues on the diagonal in decreasing order. Finally, if there are complex eigenvalues, then $A$ might be brought to the last canonical form. To this end, first normalize the real eigenvalue to zero by adding a multiple of the identity matrix, then normalize the imaginary parts of the complex eigenvalues to $\pm1$ by multiplying $A$ with a constant, make the real part of the complex eigenvalues nonnegative by possibly multiplying with $-1$, and finally pass to the real Jordan form.
\end{proof}

We shall now classify those regular convex cones $K$ whose automorphism group has one of the matrices listed in Lemma \ref{can_form} as an element $A$ of its Lie algebra $\mathfrak{aut}K$. To this end, we consider the action of the 2-parametric subgroup $G_A = \{ \lambda e^{At} \,|\, \lambda > 0,\ t \in \mathbb R \} \subset \Aut K$ on the elements $x = (x_1,x_2,x_3)^T \in \mathbb R^3$. The interior of the cone $K$ as well as its boundary have to be composed of orbits of $G_A$. Denote the canonical basis vectors of $\mathbb R^3$ by $e_1,e_2,e_3$. We consider the canonical forms in Lemma \ref{can_form} case by case.

Case 1: $e^{At} = \begin{pmatrix} 1 & t & \frac{t^2}{2} \\ 0 & 1 & t \\ 0 & 0 & 1 \end{pmatrix}$. The boundary $\partial K$ of $K$ must intersect at least one of the open half-spaces given by $x_3 > 0$, $x_3 < 0$, because $K$ has a non-empty interior. Let $x \in \partial K$ be such that $x_3 \not= 0$. The orbit of $x$, which is a subset of $\partial K$, is the conic hull of a parabola. The closure of this orbit is the boundary of a second-order cone, which must hence coincide with $K$.

Case 2: $e^{At} = \begin{pmatrix} 1 & t & 0 \\ 0 & 1 & 0 \\ 0 & 0 & 1 \end{pmatrix}$. There exists a point $x \in K$ with $x_2 \not= 0$. The orbit of this point contains a line, however, leading to a contradiction with the regularity of $K$. Hence this case is not possible.

Case 3: $e^{At} = \begin{pmatrix} e^t & te^t & 0 \\ 0 & e^t & 0 \\ 0 & 0 & 1 \end{pmatrix}$. We have $SA = AS$ for all $S \in \left\{ \left. \begin{pmatrix} \alpha & \beta & 0 \\ 0 & \alpha & 0 \\ 0 & 0 & \gamma \end{pmatrix} \right| \alpha,\beta,\gamma \in \mathbb R \right\}$. Hence $A \in \mathfrak{aut}K$ implies $A \in \mathfrak{aut}\tilde K$ with $\tilde K = S[K]$ for every invertible $S$ of this form. The boundary of $K$ must contain a point $x$ such that $x_2 \not= 0$, $x_3 \not= 0$. By applying a suitable transformation $S$ to the cone $K$, we can assume without restriction of generality that $x = (0,1,1)^T$. The orbit of $x$ is given by the set $\{\lambda(t,1,e^{-t})^T \,|\, \lambda > 0,\ t \in \mathbb R\}$, which is a 2-dimensional manifold bounded by the rays generated by $e_1$ and $e_3$. Hence $e_1,e_3 \in \partial K$, and every convex conic combination of $e_1$ and $e_3$ is an element of $K$. Let us show that such a combination cannot be in the interior of $K$. Assume the contrary. Then there exists $\tilde x \in K$ such that $\tilde x_3 > 0$, $\tilde x_2 < 0$. The orbit of $\tilde x$ is a 2-dimensional manifold bounded by the rays generated by $-e_1,e_3$. Thus both vectors $\pm e_1$ would be elements of $K$, in contradiction with the regularity of $K$. Therefore the conic convex hull of the vectors $e_1,e_3$ has to be a subset of the boundary $\partial K$. The union of this hull with the orbit of $x$ bounds a regular convex cone, which has then to coincide with $K$. It is easily seen that this cone is isomorphic to the cone obtained by the homogenization of the epigraph of the exponential function.

Case 4: $e^{At} = \diag(e^t,e^{-\mu t},e^{(\mu-1)t})$. We distinguish two subcases.

Case 4.1: All boundary points $x \in \partial K$ satisfy $x_1x_2x_3 = 0$. Then $K$ must be some orthant of $\mathbb R^3$.

Case 4.2: There exists a point $x \in \partial K$ such that $x_1x_2x_3 \not= 0$. Note that we have $SA = AS$ for all diagonal $S$. By applying a suitable invertible diagonal transformation $S$, we may assume that $x = (1,1,1)^T$. The orbit of $x$, which is given by the set $\{\lambda(e^t,e^{-\mu t},e^{(\mu-1)t})^T \,|\, \lambda > 0,\ t \in \mathbb R\}$, is then a subset of $\partial K$. We again distinguish two subcases.

Case 4.2.1: $\mu < \frac12$. Define $p = \frac{2-\mu}{1-2\mu} \in [2,\infty)$ and $q = \frac{2-\mu}{1+\mu} \in (1,2]$. Then we have $\frac1p+\frac1q = 1$. The orbit of the point $(1,\gamma,1)$, $\gamma \in \mathbb R$, is given by the equation $x_2 = \gamma x_1^{1/p}x_3^{1/q}$, with $x_1,x_3 > 0$. This is a 2-dimensional manifold bounded by the rays generated by $e_1$ and $e_3$. In particular, the orbit of $x$ is of this form, and $e_1,e_3 \in \partial K$. The hyperplanes $\spa\{e_1,e_2\}$ and $\spa\{e_2,e_3\}$ are tangent to the closure of the orbit of $x$ at $e_1$ and $e_3$, respectively, and hence must be supporting hyperplanes of $K$ at these points. This leaves three possibilities.

Case 4.2.1.1: $-e_2 \in K$. Then $K$ contains the complete orthant spanned by $e_1,-e_2,e_3$. It then must be the union of this orthant with the convex conic hull of the orbit of $x$.

Case 4.2.1.2: $-e_2 \not\in K$, but $e_1+e_3$ is in the interior of $K$. Then there exists $\alpha > 0$ such that $\tilde x = (1,-\alpha,1)^T \in \partial K$. The closure of the union of the orbits of $x$ and $\tilde x$ bounds a regular convex cone, which must coincide with $K$. If we apply the transformation $S = \diag(1,-\alpha^{-1},1)$ to $K$, then we obtain the cone defined by $\alpha^{-1}$ instead of $\alpha$. We may hence assume that $\alpha \in (0,1]$. Note that if $(\alpha,\mu) = (1,0)$, then $K$ is a second-order cone.

Case 4.2.1.3: $e_1+e_3 \in \partial K$. Then the closure of the union of the convex conic hull of $\{e_1,e_3\}$ with the orbit of $x$ bounds a regular convex cone, which must coincide with $K$.

Case 4.2.2: $\mu = -\frac12$. The the symmetry group $G_A$ contains the subgroup of transformations given by $\diag(\lambda,1,1)$, $\lambda > 0$. The cone $K$ is hence invariant with respect to multiplication of the first coordinate by an arbitrary positive number. It is not hard to see that any regular convex cone satisfying this condition must be simplicial.

Case 5: $e^{At} = \begin{pmatrix} e^{\mu t}\cos t & e^{\mu t}\sin t & 0 \\ -e^{\mu t}\sin t & e^{\mu t}\cos t & 0 \\ 0 & 0 & 1 \end{pmatrix}$. Let $x \in \partial K$ be such that $x_3 \not= 0$. Without loss of generality we may assume that $x_3 = 1$. Then the intersection of the orbit of $x$ with the plane $x_3 = 1$ is a spiral, except when $\mu = 0$, in which case it is a circle. Clearly the conic hull of a spiral cannot bound a convex cone. Thus $\mu = 0$, and $K$ is a second-order cone.

We arrive at the following classification.

{\theorem \label{th_class} Let $K \subset \mathbb R^3$ be a regular convex cone such that $\dim\Aut K \geq 2$. Then $K$ is isomorphic to exactly one of the following cones.

1. the cone obtained by the homogenization of the epigraph of the exponential function,

2. the positive orthant $\mathbb R_+^3$,

3. the cone given by $\{ x \,|\, x_2 \leq x_1^{1/p}x_3^{1/q},\ x_1 \geq 0,\ x_3 \geq 0\}$ for some $p \in [2,\infty)$, $\frac1p+\frac1q = 1$,

4. the cone given by $\{ x \,|\, -\alpha x_1^{1/p}x_3^{1/q} \leq x_2 \leq x_1^{1/p}x_3^{1/q},\ x_1 \geq 0,\ x_3 \geq 0\}$ for some $p \in [2,\infty)$, $\frac1p+\frac1q = 1$, $\alpha \in (0,1]$,

5. the cone given by $\{ x \,|\, 0 \leq x_2 \leq x_1^{1/p}x_3^{1/q},\ x_1 \geq 0,\ x_3 \geq 0\}$ for some $p \in [2,\infty)$, $\frac1p+\frac1q = 1$. }

\begin{proof}
That $K$ is isomorphic to one of the listed cones follows from the above case by case study. It rests to show that the cones in the list are mutually non-isomorphic.

Consider a compact section of the cone $K$. Its boundary is a closed curve. This curve may contain straight line segments as well as corners, i.e., discontinuities in the tangent direction to the curve. The number of line segments is 1,3,2,0,1, and the number of corners 1,3,1,0,2, respectively, for the five cases listed in the theorem. Hence no two cones belonging to different cases can be isomorphic.

Now consider the families of cones listed in cases 3 -- 5. Let $K,\tilde K$ be two cones from the same family, given by parameter values $p,\tilde p$, respectively, or $(p,\alpha),(\tilde p,\tilde\alpha)$ if the cones are from the 4th family, and let $S$ be an invertible linear map such that $S[K] = \tilde K$. The boundaries of the cones are real analytic, except at the rays generated by $e_1,e_3$, and, in the case of family 3, at the ray generated by $-e_2$. Hence $S$ maps the pair of rays generated by $e_1,e_3$ to itself, and must have one of the two forms
\[ S = \begin{pmatrix} S_{11} & S_{12} & 0 \\ 0 & S_{22} & 0 \\ 0 & S_{32} & S_{33} \end{pmatrix},\qquad S = \begin{pmatrix} 0 & S_{12} & S_{13} \\ 0 & S_{22} & 0 \\ S_{31} & S_{32} & 0 \end{pmatrix}
\]
with $S_{ij} > 0$ for $i,j = 1,3$, $S_{22} \not= 0$. The graph of the function $x_2 = x_1^{1/\tilde p}x_3^{1/\tilde q}$, which is part of the boundary of $\tilde K$, is the image of the surfaces given by $S_{22}x_2 = (S_{11}x_1+S_{12}x_2)^{1/\tilde p}(S_{32}x_2+S_{33}x_3)^{1/\tilde q}$, $S_{22}x_2 = (S_{12}x_2+S_{13}x_3)^{1/\tilde p}(S_{31}x_1+S_{32}x_2)^{1/\tilde q}$, respectively, for these two types of isomorphisms. These surfaces are part of the boundary of $K$ and must be the graphs of the function $x_2 = \gamma x_1^{1/p}x_3^{1/q}$ for $\gamma = -\alpha$ or $\gamma = 1$, however. It follows that
\begin{eqnarray*}
S_{22}\gamma x_1^{1/p}x_3^{1/q} &=& (S_{11}x_1+S_{12}\gamma x_1^{1/p}x_3^{1/q})^{1/\tilde p}(S_{32}\gamma x_1^{1/p}x_3^{1/q}+S_{33}x_3)^{1/\tilde q}, \\
S_{22}\gamma x_1^{1/p}x_3^{1/q} &=& (S_{12}\gamma x_1^{1/p}x_3^{1/q}+S_{13}x_3)^{1/\tilde p}(S_{31}x_1+S_{32}\gamma x_1^{1/p}x_3^{1/q})^{1/\tilde q},
\end{eqnarray*}
respectively, identically for $x_1,x_3 \geq 0$. This is equivalent to
\begin{eqnarray*}
S_{22}\gamma x_1^{1/p\tilde q}x_3^{1/q\tilde p} &=& (S_{11}x_1^{1/q}+S_{12}\gamma x_3^{1/q})^{1/\tilde p}(S_{32}\gamma x_1^{1/p}+S_{33}x_3^{1/p})^{1/\tilde q}, \\
S_{22}\gamma x_1^{1/p\tilde p}x_3^{1/q\tilde q} &=& (S_{12}\gamma x_1^{1/p}+S_{13}x_3^{1/p})^{1/\tilde p}(S_{31}x_1^{1/q}+S_{32}\gamma x_3^{1/q})^{1/\tilde q},
\end{eqnarray*}
respectively. Inserting $x_1 = 0$ or $x_3 = 0$, we get
\begin{eqnarray*}
S_{12}^{1/\tilde p}S_{33}^{1/\tilde q}\gamma^{1/\tilde p} = S_{11}^{1/\tilde p}S_{32}^{1/\tilde q}\gamma^{1/\tilde q} &=& 0,\\
S_{13}^{1/\tilde p}S_{32}^{1/\tilde q}\gamma^{1/\tilde q} = S_{12}^{1/\tilde p}S_{31}^{1/\tilde q}\gamma^{1/\tilde p} &=& 0,
\end{eqnarray*}
respectively. Since $\gamma,S_{11},S_{13},S_{31},S_{33} \not= 0$, we obtain $S_{12} = S_{32} = 0$ in both cases. Above conditions simplify to
\begin{eqnarray*}
S_{22}\gamma x_1^{1/p-1/\tilde p}x_3^{1/q-1/\tilde q} &=& S_{11}^{1/\tilde p}S_{33}^{1/\tilde q}, \\
S_{22}\gamma x_1^{1/p-1/\tilde q}x_3^{1/q-1/\tilde p} &=& S_{13}^{1/\tilde p}S_{31}^{1/\tilde q},
\end{eqnarray*}
respectively. Since this has to be valid for all $x_1,x_3 \geq 0$, we obtain $p = \tilde p$, $S_{22}\gamma = S_{11}^{1/\tilde p}S_{33}^{1/\tilde q}$ in the first case, and $p = \tilde q$, $S_{22}\gamma = S_{13}^{1/\tilde p}S_{31}^{1/\tilde q}$ in the second case. Note that $\tilde q \leq 2 \leq p$, and hence in the second case we have $p = \tilde q = 2$. But then we also have $p = \tilde p$. The cones in the families 3 and 5 are thus mutually non-isomorphic. Let us consider the two cases for family 4.

The isomorphism has one of the two forms
\[ S = \begin{pmatrix} S_{11} & 0 & 0 \\ 0 & \gamma^{-1}S_{11}^{1/p}S_{33}^{1/q} & 0 \\ 0 & 0 & S_{33} \end{pmatrix},\qquad S = \begin{pmatrix} 0 & 0 & S_{13} \\ 0 & \gamma^{-1}S_{13}^{1/2}S_{31}^{1/2} & 0 \\ S_{31} & 0 & 0 \end{pmatrix}.
\]
The graph of the function $x_2 = -\tilde\alpha x_1^{1/p}x_3^{1/q}$ is part of the boundary of $\tilde K$. It is the image of the surfaces given by $\gamma^{-1}S_{11}^{1/p}S_{33}^{1/q}x_2 = -\tilde\alpha(S_{11}x_1)^{1/p}(S_{33}x_3)^{1/q}$, $\gamma^{-1}S_{13}^{1/2}S_{31}^{1/2}x_2 = -\tilde\alpha(S_{13}x_3)^{1/2}(S_{31}x_1)^{1/2}$, respectively. These surfaces are part of the boundary of $K$ and must be the graphs of the function $x_2 = \tilde\gamma x_1^{1/p}x_3^{1/q}$ for $\tilde\gamma = 1$ or $\tilde\gamma = -\alpha$. Here $\gamma\tilde\gamma = -\alpha$, because if one part of the boundary of $K$ is mapped to one part of $\partial\tilde K$, then the other part of $\partial K$ is mapped to the other part of $\partial\tilde K$. Hence we obtain
\begin{eqnarray*}
\gamma^{-1}S_{11}^{1/p}S_{33}^{1/q}\tilde\gamma x_1^{1/p}x_3^{1/q} &=& -\tilde\alpha(S_{11}x_1)^{1/p}(S_{33}x_3)^{1/q}, \\
\gamma^{-1}S_{13}^{1/2}S_{31}^{1/2}\tilde\gamma x_1^{1/2}x_3^{1/2} &=& -\tilde\alpha(S_{13}x_3)^{1/2}(S_{31}x_1)^{1/2},
\end{eqnarray*}
respectively, for all $x_1,x_3 \geq 0$. This simplifies to $\gamma^{-1}\tilde\gamma = -\tilde\alpha$ in both cases. Combining with the condition $\gamma\tilde\gamma = -\alpha$, we get $\tilde\alpha = \alpha\gamma^{-2}$. If $\gamma = 1$, we obtain $\tilde\alpha = \alpha$. If $\gamma = -\alpha$, then $\tilde\alpha = \alpha^{-1}$. But $\tilde\alpha \leq 1 \leq \alpha^{-1}$, hence in this case also $\tilde\alpha = \alpha$. This proves that cones defined by different pairs $(p,\alpha)$ cannot be isomorphic.
\end{proof}

\section{Analytic solutions of affine spheres} \label{sec:affine_spheres}

In this section we compute the complete hyperbolic affine spheres which are asymptotic to the cones $K$ listed in Theorem \ref{th_class}. To this end we transform \eqref{MAequation} to an ODE by virtue of the symmetry group and solve this ODE. Then we convert this solution to a solution of the PDE \eqref{MAequation}, which finally allows us to obtain the affine sphere as an explicit analytic hypersurface immersion.

\subsection{Epigraph of the exponential function}

In this subsection we consider the cone obtained by the homogenization of the epigraph of the exponential function and compute the complete hyperbolic affine sphere which is asymptotic to the boundary of this cone. This cone is given by
\[ K = \{ (x,y,0)^T \,|\, x \leq 0,\ y \geq 0 \} \cup \left\{ (x,y,z)^T \,\left|\, \frac{y}{z} \geq \exp\frac{x}{z},\ z > 0 \right. \right\}.
\]
Invariance with respect to the unimodular subgroup of $\Aut K$ leads to the Ansatz $F(x,y,z) = -\log y-2\log z+\phi(t)$ with $t = \log\frac{y}{z}-\frac{x}{z}$ for the solution of \eqref{MAequation}. The function $\phi(t)$ must be defined for all $t > 0$. We have $\lim_{x \to -\infty} F(x,1,1) = -\infty$, $\lim_{x \to 0-} F(x,1,1) = +\infty$, because this affine half-line intersects all level surfaces of $F$. It follows that $\lim_{t \to 0} \phi(t) = +\infty$, $\lim_{t \to 0} \dot\phi(t) = -\infty$, $\lim_{t \to \infty} \phi(t) = -\infty$. The derivatives of $F$ are given by
\begin{eqnarray*}
F' &=& \begin{pmatrix} -\frac{\dot\phi}{z} \\ -\frac{1}{y}+\frac{\dot\phi}{y} \\ -\frac{2}{z}+\dot\phi\left(-\frac{1}{z}+\frac{x}{z^2}\right) \end{pmatrix},\\
F'' &=& \begin{pmatrix} \frac{\ddot\phi}{z^2} & -\frac{\ddot\phi}{yz} & \frac{\dot\phi}{z^2}-\frac{\ddot\phi}{z}\left(-\frac{1}{z}+\frac{x}{z^2}\right) \\ -\frac{\ddot\phi}{yz} & \frac{1}{y^2}-\frac{\dot\phi}{y^2}+\frac{\ddot\phi}{y^2} & \frac{\ddot\phi}{y}\left(-\frac{1}{z}+\frac{x}{z^2}\right) \\ \frac{\dot\phi}{z^2}-\frac{\ddot\phi}{z}\left(-\frac{1}{z}+\frac{x}{z^2}\right) & \frac{\ddot\phi}{y}\left(-\frac{1}{z}+\frac{x}{z^2}\right) & \frac{2}{z^2}+\dot\phi\left(\frac{1}{z^2}-\frac{2x}{z^3}\right)+\ddot\phi\left(-\frac{1}{z}+\frac{x}{z^2}\right)^2 \end{pmatrix}.
\end{eqnarray*}
Here $\dot\phi,\ddot\phi$ denote the derivatives of the function $\phi$ with respect to its argument $t$. It follows that $\det F'' = \frac{1}{y^2z^4}(2\ddot\phi-\dot\phi^2-3\dot\phi\ddot\phi+\dot\phi^3)$, and the PDE in \eqref{MAequation} becomes the ODE
\[ 2\ddot\phi-\dot\phi^2-3\dot\phi\ddot\phi+\dot\phi^3 = e^{2\phi}.
\]
It is not hard to see that the Hessian $F''$ is positive definite if and only if
\[ \begin{pmatrix} \ddot\phi & 0 & \dot\phi \\ 0 & 1-\dot\phi & \dot\phi \\ \dot\phi & \dot\phi & 2-\dot\phi \end{pmatrix} \succ 0.
\]
This happens if and only if $\dot\phi < \frac23$, $\ddot\phi > \frac{\dot\phi^2(1-\dot\phi)}{2-3\dot\phi}$. For $\dot\phi \in [0,\frac23)$ the function $\frac{\dot\phi^2(1-\dot\phi)}{2-3\dot\phi}$ is monotonically increasing. Hence if $\dot\phi(t_0) = \varepsilon$ for some $t_0$ and some $\varepsilon \in (0,\frac23)$, then for $t > t_0$ we have $\dot\phi(t) > \xi(t)$, where $\xi$ is  the solution of the equation $\dot\xi = \frac{\xi^2(1-\xi)}{2-3\xi}$ with initial condition $\xi(t_0) = \varepsilon$. This solution is given by $t - t_0 = \log\frac{1-\xi}{\xi} - \frac{2}{\xi} - \log\frac{1-\varepsilon}{\varepsilon} + \frac{2}{\varepsilon}$, and $\xi = \frac23$ is reached in finite time $t = t_0 + \log\frac12 - 3 - \log\frac{1-\varepsilon}{\varepsilon} + \frac{2}{\varepsilon}$. This contradicts the condition $\dot\phi < \frac23$. Thus we must have $\dot\phi < 0$ for all $t > 0$.

With $P = \dot\phi^2(\dot\phi-1)$ the ODE can be written as $P = e^{2\phi}+\frac{dP}{d\phi}$, which integrates to $P = -e^{2\phi}+ce^{\phi}$, with $c$ being an integration constant. It follows that $-P+\frac{c^2}{4} = (e^{\phi}-\frac{c}{2})^2$. Thus
\begin{equation} \label{expphiP}
e^{\phi} = \frac{c}{2} + \sqrt{-P+\frac{c^2}{4}},
\end{equation}
where the choice of the positive sign at the root is dictated by the conditions $P < 0$, $e^{\phi} > 0$. For $t \to \infty$ we have $e^{\phi} \to 0$, which implies that $c \leq 0$ and $P \to 0$, $\dot\phi \to 0$.

The original ODE yields
\[ \frac{dt}{d\dot\phi} = \frac{1}{\ddot\phi} = \frac{2-3\dot\phi}{e^{2\phi}+\dot\phi^2-\dot\phi^3} = \frac{2-3\dot\phi}{\left(\frac{c}{2} + \sqrt{-P+\frac{c^2}{4}}\right)^2-P} = \frac{2-3\dot\phi}{2\sqrt{-P+\frac{c^2}{4}} \left( \sqrt{-P+\frac{c^2}{4}} + \frac{c}{2} \right)}.
\]
Integrating, we obtain
\[ t = \int_{-\infty}^{\dot\phi} \frac{2-3\dot\phi}{2\sqrt{-P+\frac{c^2}{4}} \left( \sqrt{-P+\frac{c^2}{4}} + \frac{c}{2} \right)}d\dot\phi = \int_0^{\kappa} \frac{2\kappa+3}{2\sqrt{1+\kappa+\frac{c^2}{4}\kappa^3} \left( \sqrt{1+\kappa+\frac{c^2}{4}\kappa^3} + \frac{c}{2}\kappa^{3/2} \right)}d\kappa,
\]
where in the second expression for $t$ we made the substitution $\kappa = -\frac{1}{\dot\phi}$.

Let us show that the solutions for $c < 0$ lead to incomplete affine spheres. In this case the integrand above has the asymptotic expansion
\[ \frac{2\kappa+3}{2\sqrt{1+\kappa+\frac{c^2}{4}\kappa^3} \left( \sqrt{1+\kappa+\frac{c^2}{4}\kappa^3} + \frac{c}{2}\kappa^{3/2} \right)} = \frac{2\kappa+3}{\kappa+1+O(\kappa^{-1})} = 2 + \kappa^{-1} + O(\kappa^{-2})
\]
for $\kappa \to \infty$. It follows that $t = 2\kappa + \log\kappa + O(1)$ for large $\kappa$. Inverting this relation, we obtain $\kappa = \frac{t}{2} - \frac12\log t + O(1)$. On the other hand, \eqref{expphiP} yields
\[ \phi = \log\left(\frac{c}{2} + \sqrt{\frac{c^2}{4}+\kappa^{-2}-\kappa^{-3}}\right) = -\log(-c) - 2\log\kappa + O(\kappa^{-1}) = -\log(-c) - 2\log\frac{t}{2} + 2\frac{\log t}{t} + O(t^{-1}).
\]
We then get the limit
\begin{eqnarray*}
\lim_{z \to 0} F(-1,1,z) &=& \lim_{z \to 0} \left(-2\log z+\phi(z^{-1} - \log z)\right) = \lim_{z \to 0} \left(-2\log z-\log(-c) - 2\log\frac{z^{-1} - \log z}{2}\right) \\ &=& \log\left(-\frac{4}{c}\right).
\end{eqnarray*}
Thus the solution $F$ does not tend to $+\infty$ as the argument tends to the boundary of $K$.

Hence we must have $c = 0$. We then get
\[ t = \int_{-\infty}^{\dot\phi} \frac{2-3\dot\phi}{2(\dot\phi^2-\dot\phi^3)}d\dot\phi = \int_{-\infty}^{\dot\phi} \frac{1}{\dot\phi^2} - \frac{1}{2\dot\phi} + \frac{1}{2(\dot\phi-1)} d\dot\phi = -\frac{1}{\dot\phi} + \frac12\log\frac{\dot\phi-1}{\dot\phi}, \quad \phi = \frac12\log(\dot\phi^2-\dot\phi^3).
\]
Here the second relation comes from \eqref{expphiP}. Introducing the parameter $\kappa = -\frac{1}{\dot\phi}$, we obtain the parametric representation
\begin{equation} \label{implicit}
\begin{pmatrix} t \\ \phi \end{pmatrix} = \frac12\begin{pmatrix} \log(1+\kappa) + 2\kappa \\ \log(1+\kappa)-3\log\kappa \end{pmatrix}
\end{equation}
for the solution $\phi(t)$, where $\kappa$ runs through all positive reals.

The affine sphere given by the level set $F = 0$ can then be parameterized by the positive variables $\kappa,z$. Namely relations \eqref{implicit} give us
\[ -\log y-2\log z+\frac12\left(\log(1+\kappa)-3\log\kappa\right) = 0,\qquad \log\frac{y}{z}-\frac{x}{z} = \frac12\left(\log(1+\kappa) + 2\kappa\right),
\]
which yields the immersion
\[ \mathbb R_{++}^2 \ni \begin{pmatrix} \kappa \\ z \end{pmatrix} \mapsto \begin{pmatrix} x \\ y \\ z \end{pmatrix} = \begin{pmatrix} -z(3\log z + \frac32\log\kappa + \kappa) \\ z^{-2}\kappa^{-3/2}\sqrt{1+\kappa} \\ z \end{pmatrix}.
\]
Here $\mathbb R_{++}$ is the set of positive reals.

\subsection{Differential equation for the remaining cones}

The cases 2 --- 5 in Theorem \ref{th_class} can be treated in a common framework. Let $p \in [2,\infty)$, $q \in (1,2]$ be reals such that $\frac1p+\frac1q = 1$, and let $\alpha,\beta \in [0,+\infty]$. Consider the open cone
\[ K^o = \left\{ (x,y,z)^T \,|\, -\alpha x^{1/p}y^{1/q} < z < \beta x^{1/p}y^{1/q},\ x > 0,\ y > 0 \right\}
\]
and let $K$ be its closure. Then $K$ is a regular convex cone, except for the cases $(\alpha,\beta) = (0,0)$ and $(\alpha,\beta) = (+\infty,+\infty)$ which we henceforth exclude. The cases 2 --- 5 in Theorem \ref{th_class} are given by the values $(\alpha,\beta) = (0,\infty),(\infty,1),(\alpha,1),(0,1)$, respectively, under the identification $(x,y,z) \leftrightarrow (x_1,x_2,x_3)$ of coordinates.

Invariance with respect to the unimodular subgroup of $\Aut K$ of the solution of \eqref{MAequation} leads to the Ansatz
\begin{equation} \label{ansatz}
F(x,y,z) = -\frac{p+1}{p}\log x - \frac{q+1}{q}\log y + \phi(x^{-1/p}y^{-1/q}z),
\end{equation}
where $\phi: (-\alpha,\beta) \to \mathbb R$ is a function of a scalar variable $t = x^{-1/p}y^{-1/q}z$.

Knowledge of the solution $\phi(t)$ permits to recover the affine sphere given by the level surface $F = 0$, as follows. Suppose this solution is given parametrically by functions $t(\xi),\phi(\xi)$, with $\xi$ a real parameter taking values in some interval $\Xi \subset \mathbb R$. Introduce the variable $\mu = \log x - \log y \in \mathbb R$. Then the relations $\frac{p+1}{p}\log x + \frac{q+1}{q}\log y = \phi$, $z = x^{1/p}y^{1/q}t$ yield $\log x = \frac{\phi}{3} + \frac{q+1}{3q}\mu$, $\log y = \frac{\phi}{3} - \frac{p+1}{3p}\mu$, $z = e^{\frac{\phi}{3}}e^{-\frac{\mu(p-q)}{3(p+q)}}t$. The affine sphere is then given by the immersion
\begin{equation} \label{general_immersion}
\Xi \times \mathbb R \ni \begin{pmatrix} \xi \\ \mu \end{pmatrix} \mapsto \begin{pmatrix} x \\ y \\ z \end{pmatrix} = e^{\frac{\phi(\xi)}{3}} \begin{pmatrix} e^{\frac{q+1}{3q}\mu} \\ e^{-\frac{p+1}{3p}\mu} \\ e^{-\frac{\mu(p-q)}{3(p+q)}}t(\xi) \end{pmatrix}.
\end{equation}

The derivatives of $F$ are given by
\begin{eqnarray*}
F' &=& \begin{pmatrix} -\frac{p+1+t\dot\phi}{px} \\ -\frac{q+1+t\dot\phi}{qy} \\ \frac{t\dot\phi}{z} \end{pmatrix}, \\
F'' &=& \begin{pmatrix} \frac{(p+1)(p+t\dot\phi)+t^2\ddot\phi}{p^2x^2} & \frac{t\dot\phi+t^2\ddot\phi}{pqxy} & -\frac{t\dot\phi+t^2\ddot\phi}{pxz} \\
\frac{t\dot\phi+t^2\ddot\phi}{pqxy} & \frac{(q+1)(q+t\dot\phi)+t^2\ddot\phi}{q^2y^2} & -\frac{t\dot\phi+t^2\ddot\phi}{qyz} \\
-\frac{t\dot\phi+t^2\ddot\phi}{pxz} & -\frac{t\dot\phi+t^2\ddot\phi}{qyz} & \frac{t^2\ddot\phi}{z^2} \end{pmatrix}.
\end{eqnarray*}
Here $\dot\phi,\ddot\phi$ denote the derivatives of the function $\phi$ with respect to its argument $t$. It is not hard to see that the Hessian $F''$ is positive definite if and only if
\[ \begin{pmatrix} (p+1)p+(p-1)t\dot\phi & -t\dot\phi & \dot\phi \\
-t\dot\phi & (q+1)q+(q-1)t\dot\phi & \dot\phi \\
\dot\phi & \dot\phi & \ddot\phi \end{pmatrix} \succ 0.
\]
This happens if and only if
\begin{equation} \label{convexity_condition}
t\dot\phi > -\frac23(p+q)-\frac13,\qquad \ddot\phi > \frac{\dot\phi^2(p+q-1+t\dot\phi)}{2(p+q)+1+3t\dot\phi}.
\end{equation}

Using $pq = p+q$, we obtain that $\det F'' = \frac{t^2}{pqx^2y^2z^2}\left( \ddot\phi(2(p+q)+1+3t\dot\phi) - \dot\phi^2(p+q-1+t\dot\phi) \right)$. The PDE in \eqref{MAequation} becomes the ODE
\[ \ddot\phi(2(p+q)+1+3t\dot\phi) - \dot\phi^2(p+q-1+t\dot\phi) = (p+q)e^{2\phi}.
\]

Now we have
\begin{eqnarray*}
\frac{d}{dt}\left(e^{-\phi}\dot\phi(t\dot\phi+p+1)(t\dot\phi+q+1)\right) &=& e^{-\phi}\left( \ddot\phi(2(p+q)+1+3t\dot\phi) - \dot\phi^2(p+q-1+t\dot\phi) \right)(1+t\dot\phi), \\
\frac{d}{dt}\left((p+q)te^{\phi}\right) &=& (p+q)e^{\phi}(1+t\dot\phi).
\end{eqnarray*}
It follows that
\begin{equation} \label{main_ODE}
e^{-\phi}\dot\phi(t\dot\phi+p+1)(t\dot\phi+q+1) = (p+q)te^{\phi} + c,
\end{equation}
where $c$ is an integration constant.

We shall first investigate the solution $\phi(t)$ for $t > 0$. Let us introduce the variables $\tau = \log t$, $\varphi = \phi + \tau$, $\xi = \frac{d\varphi}{d\tau} = 1+t\dot\phi$. We then have $e^{\varphi} = te^{\phi}$, and the equation becomes $(\xi-1)(\xi+p)(\xi+q) = e^{\varphi}((p+q)e^{\varphi} + c)$. Writing shorthand $P = (\xi-1)(\xi+p)(\xi+q) = \xi^3 + \xi^2(p+q-1) - p - q$, we obtain
\begin{equation} \label{expphi_exp}
e^{\varphi} = \frac{-c+\sigma\sqrt{c^2+4(p+q)P}}{2(p+q)}
\end{equation}
with $\sigma = \pm1$. Differentiating both sides with respect to $\tau$, we obtain $\xi e^{\varphi} = \frac{d\xi}{d\tau}\frac{dP}{d\xi}\frac{\sigma}{\sqrt{c^2+4(p+q)P}}$,
which after substitution of $e^{\varphi}$ by virtue of \eqref{expphi_exp} and resolution with respect to $\frac{d\xi}{d\tau}$ yields
\begin{equation} \label{xi_tau_ODE}
\frac{d\xi}{d\tau} = \frac{-c\sigma\sqrt{c^2+4(p+q)P}+c^2+4(p+q)P}{2(p+q)(3\xi + 2(p+q-1))}.
\end{equation}
The right-hand side of this equation depends on $\xi$ only. Since $\xi$ is an analytic function of $\tau$, we obtain two cases. Either $\xi$ is constant, or $\frac{d\xi}{d\tau} \not= 0$ almost everywhere.

Let us first consider the case $\xi \equiv c_1$. Then $\dot\phi = \frac{c_1-1}{t}$, $\phi = (c_1-1)\log t + c_2$, and $F(x,y,z) = -\frac{p+c_1}{p}\log x - \frac{q+c_1}{q}\log y + (c_1-1)\log z + c_2$. We have $F'' \succ 0$ if and only if $-q < c_1 < 1$, and $\det F'' = e^{2F}$ if and only if $c_1 = c_2 = 0$. Then $K$ is the positive orthant, $\beta = +\infty$, and $\alpha = 0$. The function $F$ is the well-known solution $-\log(xyz)$, and the affine spheres asymptotic to the orthant are the level surfaces of the product $xyz$.

\subsection{Solution of the differential equation}

In this subsection we solve ODE \eqref{xi_tau_ODE} in the case $\xi \not\equiv const$. We can then use $\xi$ as an independent variable, which yields
\begin{eqnarray*}
\frac{d\tau}{d\xi} &=& \frac{2(p+q)(3\xi + 2(p+q-1))}{\sqrt{c^2+4(p+q)P}(-c\sigma+\sqrt{c^2+4(p+q)P})} \\ &=& \frac{(3\xi + 2(p+q-1))(\sqrt{c^2+4(p+q)P}+c\sigma)}{2P\sqrt{c^2+4(p+q)P}} \\ &=& \left( \frac{1}{\xi-1} - \frac{1}{p(\xi+p)} - \frac{1}{q(\xi+q)} \right)\left(\frac12+\frac{c\sigma}{2\sqrt{c^2+4(p+q)P}}\right).
\end{eqnarray*}

We have $\frac{d\xi}{d\tau} = t^2\ddot\phi + t\dot\phi$. From \eqref{convexity_condition} we then get $\xi > -\frac23(p+q-1)$. Note that $\xi = -\frac23(p+q-1)$ is the location of the local maximum of the cubic polynomial $P$. At this maximum we have $P = \frac{(p+q-4)(2(p+q)+1)^2}{27}$. The local minimum of $P$ is located at $\xi = 0$, with value $P = -(p+q)$.

By continuity of $\varphi$ the sign $\sigma$ can change only if $c^2+4(p+q)P = 0$. The cubic polynomial $c^2+4(p+q)P$ in $\xi$ always has a real root $\xi_3$ which does not exceed the smallest real root of $P$, i.e., $\xi_3 \leq -p$. The other two roots are real if and only if $|c| \leq 2(p+q)$. If this is the case, then we denote the larger of the two roots by $\xi_1$ and the smaller one by $\xi_2$. These roots are located between the two larger real roots of $P$ and are separated by the local minimum of $P$, i.e., $-q \leq \xi_2 \leq 0 \leq \xi_1 \leq 1$. If $|c| = 2(p+q)$, then $\xi_1 = \xi_2 = 0$ is a double root, and $\xi_3 = -(p+q-1)$.

By virtue of \eqref{expphi_exp} we must have $c^2+4(p+q)P \geq 0$ and $\sigma\sqrt{c^2+4(p+q)P} > c$. If $|c| > 2(p+q)$, then the first condition is satisfied automatically. In the case $|c| \leq 2(p+q)$ it is equivalent to $-\frac23(p+q-1) < \xi \leq \xi_2$ or $\xi_1 \leq \xi$. The second condition implies that if $c \geq 0$, then $P > 0$ and $\sigma = 1$. Moreover, $P \geq 0$ implies $\sigma = 1$. It also follows that $\frac{d\xi}{d\tau} \geq 0$ if $\sigma = 1$ and $\frac{d\xi}{d\tau} \leq 0$ if $\sigma = -1$.

\medskip

We now investigate the qualitative behaviour of the solutions. The derivative $\frac{d\tau}{d\xi}$ becomes zero at $\xi = -\frac23(p+q-1)$ and $\xi = +\infty$. It may become infinite at the roots of the polynomial $P$ or at the roots of $c^2+4(p+q)P$. Note that the ODE is invariant with respect to shifts of the variable $\tau$ by an additive constant. We shall now analyze the behaviour of the solution in the neighbourhood of the above-mentioned points.

$\bm{\xi \to +\infty}$. In this case $P > 0$, and hence $\sigma = 1$. We have $\frac{d\tau}{d\xi} = \frac32\xi^{-2} + O(\xi^{-3})$, $\tau = -\frac32\xi^{-1}+ O(\xi^{-2}) + const$, and $\xi$ escapes to infinity in finite time $\tau = \tau^*$. Moreover, we have $\xi = \frac32\frac{1}{\tau^*-\tau} + O(1)$, $\dot\phi = \frac32\frac{1}{t^*-t} + O(1)$, and $\phi = -\frac32\log(t^*-t) + O(1)$, where $t^* = e^{\tau^*}$. Hence the solution $\phi(t)$ also escapes to $+\infty$ in finite time.

$\bm{\xi \to 1^+}$. For $\xi > 1$ we have $P > 0$ and hence $\sigma = 1$. If $c > 0$, then $\frac{d\tau}{d\xi} = \frac{1}{\xi-1} + O(1)$, $\tau = \log(\xi-1) + O(1)$, and $\tau$ escapes to $-\infty$. Moreover, $\log\dot\phi = O(1)$, and $\lim_{t \to 0}\dot\phi$ is finite and positive. If $c = 0$, then $\frac{d\tau}{d\xi} = \frac{1}{2(\xi-1)} + O(1)$, $\tau = \frac12\log(\xi-1) + O(1)$, and $\tau$ also escapes to $-\infty$. In this case we have $\log t = \log\dot\phi + O(1)$, and $\lim_{t \to 0}\dot\phi = 0$. If $c < 0$, then $\lim_{\xi\to1^+}\frac{d\tau}{d\xi} = \frac{(p+q)(p+1)(q+1)}{c^2}$ is finite and positive, and the solution continues beyond $\xi = 1$.

$\bm{\xi \to 1^-}$. In this case $P < 0$ and hence $c < 0$. If $\sigma = 1$, then $\lim_{\xi\to1^-}\frac{d\tau}{d\xi} = \frac{(p+q)(p+1)(q+1)}{c^2}$ is finite and positive, and the solution continues beyond $\xi = 1$. If $\sigma = -1$, then $\frac{d\tau}{d\xi} = \frac{1}{\xi-1} + O(1)$, $\tau = \log(1-\xi) + O(1)$, and $\tau$ escapes to $-\infty$. In this case $\log(-\dot\phi) = O(1)$, and $\lim_{t \to 0}\dot\phi$ is finite and negative.

$\bm{\xi \to \xi_1^+}$. Then $|c| \leq 2(p+q)$. We assume $c \not= 0$, otherwise $\xi_1 = 1$. We then have $P < 0$ and hence $c < 0$. It follows that $\frac{d\tau}{d\xi} = \frac{2(p+q)(3\xi_1 + 2(p+q-1))}{-c\sigma\sqrt{c^2+4(p+q)P}} + O(1)$. If $c > -2(p+q)$, then $\xi_1$ is a single root of the polynomial $c^2+4(p+q)P$, $\frac{d\tau}{d\xi} = O((\xi-\xi_1)^{-1/2})$, $\tau = O(\sqrt{\xi-\xi_1}) + const$, and $\xi = \xi_1$ is reached in finite time $\tau = \tau_1$. Moreover, $\xi = \xi_1 + O((\tau-\tau_1)^2)$, and the passage through the point $\tau = \tau_1$ is accompanied by a switch in the sign $\sigma$. If $c = -2(p+q)$, then $\xi_1 = \xi_2 = 0$ is a double root of $c^2+4(p+q)P$. In this case $\frac{d\tau}{d\xi} = \frac{\sqrt{p+q-1}}{\sigma\xi\sqrt{p+q}} + O(1)$, $\tau = \frac{\sqrt{p+q-1}}{\sigma\sqrt{p+q}}\log\xi + O(1)$, and $\tau$ escapes to $\mp\infty$ for $\sigma = \pm1$. We get $\log(1+t\dot\phi) = \sigma\sqrt{\frac{p+q}{p+q-1}}\log t + O(1)$, $t\dot\phi = -1+O(t^{\sigma\sqrt{\frac{p+q}{p+q-1}}})$, and hence $\phi = -\log t + O(1)$.

$\bm{\xi \to \xi_2^-}$. Then $|c| \leq 2(p+q)$. We also assume $c \not= 0$, otherwise $\xi_2 = -q$, which will be considered below. Exactly as in the previous case, we have $c < 0$. If $c > -2(p+q)$, then $\xi = \xi_2$ is reached in finite time $\tau = \tau_2$, and the passage through this point is accompanied by a switch in the sign $\sigma$. If $c = -2(p+q)$, then $\xi_1 = \xi_2 = 0$ and $\tau$ escapes to $\pm\infty$ for $\sigma = \pm1$. As in the case $\xi \to \xi_1^+$, we have $\phi = -\log t + O(1)$.

$\bm{\xi \to -q^+}$. As in the case $\xi\to1^-$, we have $c < 0$. If $\sigma = 1$, then $\lim_{\xi\to -q^+}\frac{d\tau}{d\xi} = \frac{(p^2-q^2)(q+1)}{qc^2}$ is finite. If $p > q$, then the limit is positive and the solution continues beyond $\xi = -q$. If $p = q$, then $\lim_{\xi\to -q^+}\frac{d\tau}{d\xi} = 0$, and $\xi = -q = -\frac23(p+q-1)$ is reached in finite time $\tau$. If $\sigma = -1$, then $\frac{d\tau}{d\xi} = -\frac{1}{q(\xi+q)} + O(1)$, $\tau = -\frac1q\log(\xi+q) + O(1)$, $t\dot\phi = -1-q+O(t^{-q})$ for $p > q$, and $\frac{d\tau}{d\xi} = -\frac{1}{\xi+2} + O(1)$, $\tau = -\log(\xi+2) + O(1)$, $t\dot\phi = -3+O(t^{-1})$ for $p = q$. In either case, $\tau$ escapes to $+\infty$ and $\phi = -(q+1)\log t + O(1)$.

$\bm{\xi \to -q^-}$. The condition $\xi > -\frac23(p+q-1)$ implies that $-\frac23(p+q-1) < -q$, i.e., $p > q$. As in the case $\xi\to1^+$, we have $\sigma = 1$. If $c > 0$, then $\frac{d\tau}{d\xi} = -\frac{1}{q(\xi+q)} + O(1)$, $\tau = -\frac1q\log(-\xi-q) + O(1)$, and $\tau$ escapes to $+\infty$. Moreover, we have $t\dot\phi = -1-q+O(t^{-q})$, and $\phi = -(q+1)\log t + O(1)$. If $c = 0$, then $\frac{d\tau}{d\xi} = -\frac{1}{2q(\xi+q)} + O(1)$, $\tau = -\frac{1}{2q}\log(-\xi-q) + O(1)$, and $\tau$ also escapes to $+\infty$. We have $t\dot\phi = -1-q+O(t^{-2q})$, and $\phi = -(q+1)\log t + O(1)$ too. If $c < 0$, then $\lim_{\xi\to -q^-}\frac{d\tau}{d\xi} = \frac{(p^2-q^2)(q+1)}{qc^2}$ is finite and positive, and the solution continues beyond $\xi = -q$.

$\bm{\xi \to -\frac23(p+q-1)}$. If $p = q$, then $-\frac23(p+q-1) = -q$, which case was considered above. If $p > q$, then we have $P > 0$ and hence $\sigma = 1$. The derivative $\frac{d\tau}{d\xi}$ tends to zero asymptotically proportionally to $\xi + \frac23(p+q-1)$. It follows that $\tau = const + O(\sqrt{\xi + \frac23(p+q-1)})$, and $\xi = -\frac23(p+q-1)$ is reached in finite time $\tau$.

\medskip

We are now able to construct the solutions of \eqref{xi_tau_ODE} for different values of $c$.

$\bm{c \geq 0}$. In this case $\sigma = 1$. There exists a strictly monotonically increasing solution $\xi(\tau)$ which is defined on the semi-infinite interval $(-\infty,\tau^*)$, $\tau^* \in \mathbb R$ arbitrary. We have $\lim_{\tau \to -\infty}\xi = 1$, $\lim_{t \to 0}\dot\phi \geq 0$, with equality holding if and only if $c = 0$, and $\lim_{\tau\to\tau^*}\xi = +\infty$, $\lim_{t \to t^*}\phi = +\infty$. The solution is given by
\begin{equation} \label{tau_positive_int}
\tau = \tau^* - \int_{\xi}^{+\infty} \left( \frac{1}{\xi-1} - \frac{1}{p(\xi+p)} - \frac{1}{q(\xi+q)} \right)\left(\frac12+\frac{c}{2\sqrt{c^2+4(p+q)P}}\right)d\xi.
\end{equation}
If $p > q$, then there exists another strictly monotonically increasing solution, defined on a semi-infinite interval $(\tau^*,+\infty)$, with $\lim_{\tau \to \tau^*}\xi = -\frac23(p+q-1)$ and $\lim_{\tau\to+\infty}\xi = -q$. This second solution reaches $\xi = -\frac23(p+q-1)$ in finite time and hence violates \eqref{convexity_condition}.

$\bm{-2(p+q) < c < 0}$. For $\sigma = 1$ we have a strictly monotonically increasing solution $\xi(\tau)$ escaping to $+\infty$ at finite time $\tau = \tau^*$. This solution is given by \eqref{tau_positive_int} and obeys $\lim_{t \to t^*}\phi = +\infty$. However, contrary to the case $c \geq 0$ considered above, this solution can be continued beyond $\xi = 1$ up to the root $\xi = \xi_1$, which is reached in finite time
\[ \tau_1 = \tau^* - \int_{\xi_1}^{+\infty} \left( \frac{1}{\xi-1} - \frac{1}{p(\xi+p)} - \frac{1}{q(\xi+q)} \right)\left(\frac12+\frac{c}{2\sqrt{c^2+4(p+q)P}}\right)d\xi.
\]
At $\tau = \tau_1$ the sign $\sigma$ switches to $-1$ and for $\tau < \tau_1$ the solution $\xi(\tau)$ becomes strictly monotonically decreasing, reaching $\xi = 1$ from below as $\tau \to -\infty$. Moreover, we have $\lim_{t \to 0}\dot\phi < 0$. This second branch of the solution is given by
\begin{equation} \label{tau_negative_int}
\tau = \tau_1 + \int_{\xi_1}^{\xi} \left( \frac{1}{\xi-1} - \frac{1}{p(\xi+p)} - \frac{1}{q(\xi+q)} \right)\left(\frac12-\frac{c}{2\sqrt{c^2+4(p+q)P}}\right)d\xi.
\end{equation}
There exists also another solution $\xi(\tau)$. On the interval $(\tau_2,+\infty)$ it is strictly monotonically decreasing from $\xi = \xi_2$ to $\xi = -q$ with $\sigma = -1$. At $\tau = \tau_2$ the sign $\sigma$ switches to $+1$, and for $\tau < \tau_2$ the solution is strictly monotonically increasing. However, it reaches $\xi = -\frac23(p+q-1)$ in finite time, and hence violates \eqref{convexity_condition}.

$\bm{c = -2(p+q)}$. In this case the strictly monotonically increasing solution $\xi(\tau)$ given by \eqref{tau_positive_int} is valid on the whole semi-infinite interval $(-\infty,\tau^*)$, and obeys $\lim_{\tau\to-\infty}\xi = 0$. It also satisfies $\lim_{t \to t^*}\phi = +\infty$, $\phi(t) = -\log t + O(1)$ for $t \to 0$. There exist two strictly monotonically decreasing solutions with $\sigma = -1$ which are defined for all $\tau \in \mathbb R$. For one of them $\xi$ runs through the interval $(0,1)$, for the other one through the interval $(-q,0)$. For $\xi \to 0$ these solutions obey $\phi(t) = -\log t + O(1)$. For $\xi \to 1$ the derivative $\dot\phi$ tends to a finite negative value. For $\xi \to -q$ we have $\phi = -(q+1)\log t + O(1)$. There exists also a second strictly monotonically increasing solution with $\sigma = 1$, which is defined on a semi-infinite interval $(\tau^*,+\infty)$. For this solution, $\lim_{\tau \to \tau^*}\xi = -\frac23(p+q-1)$ and $\lim_{\tau\to+\infty}\xi = 0$. Hence it reaches $\xi = -\frac23(p+q-1)$ in finite time and violates \eqref{convexity_condition}.

$\bm{c < -2(p+q)}$. In this case the strictly monotonically increasing solution $\xi(\tau)$ given by \eqref{tau_positive_int} extends up to $\xi = -\frac23(p+q-1)$, which is reached in finite time. Hence this solution violates \eqref{convexity_condition}. There exists also a strictly monotonically decreasing solution with $\sigma = -1$, which is defined for all $\tau \in \mathbb R$. On this solution $\xi$ runs through the interval $(-q,1)$. The limit $\lim_{t\to0}\dot\phi$ is finite and negative, and for $t\to+\infty$ we have $\phi = -(q+1)\log t + O(1)$.

We shall now consider the solutions with asymptotic $\phi(t) = -(q+1)\log t + O(1)$ for $t \to +\infty$, which exist for $c \leq -2(p+q)$. In this case $\beta = +\infty$, and the positive orthant is contained in the cone $K$. In particular, $(1,y,1)^T \in K$ for all $y \geq 0$, and $(1,0,1)^T \in \partial K$. It follows that $\lim_{y \to 0^+}F(1,y,1) = +\infty$. For $x = z = 1$ we have $t = x^{-1/p}y^{-1/q}z = y^{-1/q}$, and hence $y = t^{-q}$. This yields $F(1,y,1) = -\frac{q+1}{q}\log y + \phi(t) = (q+1)\log t + \phi(t) = O(1)$ for $t \to +\infty$, leading to a contradiction. Thus the solution with asymptotic $\phi(t) = -(q+1)\log t + O(1)$ do not yield complete hyperbolic affine spheres. In particular, we cannot have $c < -2(p+q)$.

\subsection{Construction of the affine spheres}

We are now in a position to construct the solutions $\phi(t)$ on the whole interval $(-\alpha,\beta)$, and from these the complete hyperbolic affine spheres which are asymptotic to the boundary of the cones listed in Theorem \ref{th_class}. Equation \eqref{main_ODE} possesses the one-parametric symmetry group $t \mapsto \lambda t$, $\phi \mapsto \phi - \log\lambda$, $\lambda > 0$. If $\phi(t)$ is a solution of the equation for the value $c$ of the integration constant, then $\tilde\phi(t) = \phi(-t)$ is a solution for the value $\tilde c = -c$. We may thus assume without restriction of generality that $c \leq 0$ and either $\beta = 1$ or $\beta = +\infty$. We again consider several cases.

$\bm{c = 0}$. Then $\beta = t^* = 1$, and \eqref{tau_positive_int} becomes
\[ \log t = -\frac12\int_{\xi}^{+\infty} \left( \frac{1}{\xi-1} - \frac{1}{p(\xi+p)} - \frac{1}{q(\xi+q)} \right)d\xi = \frac{1}{2p}\log\frac{\xi-1}{\xi+p}+\frac{1}{2q}\log\frac{\xi-1}{\xi+q}.
\]
Moreover, \eqref{expphi_exp} yields $e^{\varphi} = te^{\phi} = \sqrt{\frac{P}{p+q}}$, which leads to a parametric representation of the solution $\phi(t)$ for $t \in [0,1)$. Namely, we have
\[ \begin{pmatrix} t \\ \phi \end{pmatrix} = \begin{pmatrix} (\xi+p)^{-\frac{1}{2p}}(\xi+q)^{-\frac{1}{2q}}(\xi-1)^{\frac12} \\ -\frac12\log(p+q)+\frac{p+1}{2p}\log(\xi+p)+\frac{q+1}{2q}\log(\xi+q) \end{pmatrix},
\]
with the parameter $\xi$ running from 1 to $+\infty$.

Note that $t^2 = (\xi+p)^{-\frac1p}(\xi+q)^{-\frac1q}(\xi-1)$ is analytic in $\xi$ in the neighbourhood of $\xi = 1$, and $\left.\frac{dt^2}{d\xi}\right|_{\xi = 1} = (p+1)^{-\frac1p}(q+1)^{-\frac1q} \not= 0$. Hence $\xi$ is an analytic function of $\tilde t = t^2$ in the neighbourhood of $\tilde t = 0$. But then $\phi$, which is analytic in $t$, must be an even function of $t$. In particular, $\alpha = \beta = 1$, and for negative $t$ the solution $\phi(t)$ is given parametrically by
\[ \begin{pmatrix} t \\ \phi \end{pmatrix} = \begin{pmatrix} -(\xi+p)^{-\frac{1}{2p}}(\xi+q)^{-\frac{1}{2q}}(\xi-1)^{\frac12} \\ -\frac12\log(p+q)+\frac{p+1}{2p}\log(\xi+p)+\frac{q+1}{2q}\log(\xi+q) \end{pmatrix},
\]
where $\xi \in (1,+\infty)$. Introducing a parameter $\zeta \in (-\infty,+\infty)$ such that $\xi = 1+\zeta^2$, we can represent the solution $\phi(t)$ on the whole interval $t \in (-1,1)$ by
\[ \begin{pmatrix} t \\ \phi \end{pmatrix} = \begin{pmatrix} (\zeta^2+p+1)^{-\frac{1}{2p}}(\zeta^2+q+1)^{-\frac{1}{2q}}\zeta \\ -\frac12\log(p+q)+\frac{p+1}{2p}\log(\zeta^2+p+1)+\frac{q+1}{2q}\log(\zeta^2+q+1) \end{pmatrix}.
\]

By \eqref{general_immersion} the affine sphere given by the level surface $F = 0$ can be represented as the hypersurface immersion
\[ \mathbb R^2 \ni \begin{pmatrix} \zeta \\ \mu \end{pmatrix} \mapsto \begin{pmatrix} x \\ y \\ z \end{pmatrix} = (p+q)^{-\frac16}\begin{pmatrix} (\zeta^2+p+1)^{\frac{p+1}{6p}}(\zeta^2+q+1)^{\frac{q+1}{6q}}\exp\left(\frac{q+1}{3q}\mu\right) \\ (\zeta^2+p+1)^{\frac{p+1}{6p}}(\zeta^2+q+1)^{\frac{q+1}{6q}}\exp\left(-\frac{p+1}{3p}\mu\right) \\ (\zeta^2+p+1)^{\frac{p-2}{6p}}(\zeta^2+q+1)^{\frac{q-2}{6q}}\zeta\exp\left(-\frac{p-q}{3(p+q)}\mu\right) \end{pmatrix}.
\]
This affine sphere is asymptotic to the cone given in case 4 of Theorem \ref{th_class}, with $\alpha = 1$. This cone is well-known from conic optimization, where it is called {\it power cone}. The canonical potentials for a family of cones including the power cone have been already computed in \cite{Hildebrand12b}.

$\bm{-2(p+q) < c < 0}$. In this case also $\beta = t^* = 1$, and \eqref{tau_positive_int} becomes
\[ \log t = -\int_{\xi}^{+\infty} \left( \frac{1}{\xi-1} - \frac{1}{p(\xi+p)} - \frac{1}{q(\xi+q)} \right)\left(\frac12+\frac{c}{2\sqrt{c^2+4(p+q)P}}\right)d\xi.
\]
This solution is valid for $\xi \in [\xi_1,+\infty)$. If we introduce the variable $\zeta$ by $\xi = \zeta^2 + \xi_1$, then we can combine this solution with the second branch \eqref{tau_negative_int}, yielding
\begin{equation} \label{tau_int_middle}
\log t = -\int_{\zeta}^{+\infty} G(\zeta) d\zeta,
\end{equation}
where
\begin{equation} \label{G_expr}
G = \left( \frac{1}{\zeta^2 + \xi_1 - 1} - \frac{1}{p(\zeta^2 + \xi_1 + p)} - \frac{1}{q(\zeta^2 + \xi_1 + q)} \right) \left(\zeta+\frac{c}{2\sqrt{(p+q)(\zeta^2 + \xi_1 - \xi_2)(\zeta^2 + \xi_1 - \xi_3)}}\right).
\end{equation}
Here we used that $c^2+4(p+q)P = 4(p+q)(\xi-\xi_1)(\xi-\xi_2)(\xi-\xi_3)$. Solution \eqref{tau_int_middle} is valid for $\zeta \in (-\sqrt{1-\xi_1},+\infty)$. On this interval the integrand \eqref{G_expr} is a positive analytic function. At $\zeta = -\sqrt{1-\xi_1}$ it has a simple pole with residue 1. Finally, by virtue of \eqref{tau_positive_int} the solution for negative $t$ is given by
\begin{eqnarray} \label{tau_int_left}
\log(-t) &=& \log\alpha - \int_{\xi}^{+\infty} \left( \frac{1}{\xi-1} - \frac{1}{p(\xi+p)} - \frac{1}{q(\xi+q)} \right)\left(\frac12-\frac{c}{2\sqrt{c^2+4(p+q)P}}\right)d\xi \nonumber\\ &=& \log\alpha + \int_{-\infty}^{\zeta} G(\zeta) d\zeta,
\end{eqnarray}
which is valid for $\xi \in (1,+\infty)$ and $\zeta \in (-\infty,-\sqrt{1-\xi_1})$, respectively. On this interval $G(\zeta)$ is a negative analytic function.

The constant $\alpha > 0$ must be chosen such that the limits $\lim_{t \to 0^{\pm}}\dot\phi$ coincide. We have $\sign t = \sign(\zeta+\sqrt{1-\xi_1})$, and hence
\begin{eqnarray*} 
\lim_{t \to 0}\dot\phi &=& \lim_{\xi\to1}\frac{\xi-1}{t} = \lim_{\zeta\to-\sqrt{1-\xi_1}}\frac{(\zeta-\sqrt{1-\xi_1})|\zeta+\sqrt{1-\xi_1}|}{e^{\tau}} \\ &=& -2\sqrt{1-\xi_1}\exp\left(\lim_{\zeta\to-\sqrt{1-\xi_1}}\log|\zeta+\sqrt{1-\xi_1}|-\tau\right).
\end{eqnarray*}
It follows that 
\[ \lim_{\varepsilon\to0^+} \left( \log\alpha + \int_{-\infty}^{-\sqrt{1-\xi_1}-\varepsilon} G(\zeta)d\zeta - \log\varepsilon \right) = \lim_{\varepsilon\to0^+} \left( -\int_{-\sqrt{1-\xi_1}+\varepsilon}^{+\infty} G(\zeta)d\zeta - \log\varepsilon \right).
\]
This finally yields
\begin{equation} \label{alpha_expr}
\log\alpha = -\int_0^{+\infty} G(\zeta-\sqrt{1-\xi_1}) + G(-\zeta-\sqrt{1-\xi_1}) d\zeta.
\end{equation}
Now the pole in the function $G$ cancels out and the integral on the right-hand side is finite. It is not difficult to check that the integrand is positive and hence $\alpha < 1$. The integrals can be expressed by elliptic functions \cite{ByrdFriedman54}, see also \cite{AbramowitzStegun}. 

Equation \eqref{expphi_exp} becomes
\begin{equation} \label{phi_param}
\phi = \log\frac{-c+2\zeta\sqrt{(p+q)(\zeta^2 + \xi_1 - \xi_2)(\zeta^2 + \xi_1 - \xi_3)}}{2t(p+q)}.
\end{equation}
Together with \eqref{phi_param} equations \eqref{tau_int_middle},\eqref{tau_int_left} yield a parametric representation of the solution $\phi(t)$, with $G$ given by \eqref{G_expr}, $\alpha$ given by \eqref{alpha_expr}, and $t = 0$ for $\zeta = -\sqrt{1-\xi_1}$. The parameter $\zeta$ runs through $\mathbb R$. By virtue of \eqref{general_immersion} this representation allows to construct the affine sphere given by the level surface $F = 0$ by the immersion
\[ \mathbb R^2 \ni \begin{pmatrix} \zeta \\ \mu \end{pmatrix} \mapsto \begin{pmatrix} x \\ y \\ z \end{pmatrix} = e^{\frac{\phi(\zeta)}{3}} \begin{pmatrix} e^{\frac{q+1}{3q}\mu} \\ e^{-\frac{p+1}{3p}\mu} \\ e^{-\frac{\mu(p-q)}{3(p+q)}}t(\zeta) \end{pmatrix}.
\]
This affine sphere is asymptotic to the cone given in case 4 of Theorem \ref{th_class}, for $\alpha \not= 1$.

$\bm{c = -2(p+q),\beta = 1}$. We have $t^* = \beta = 1$ and $\alpha = 0$. Equation \eqref{tau_positive_int} becomes
\begin{eqnarray*}
\log t &=& -\int_{\xi}^{+\infty} \left( \frac{1}{\xi-1} - \frac{1}{p(\xi+p)} - \frac{1}{q(\xi+q)} \right)\left(\frac12-\frac{\sqrt{p+q}}{2\xi\sqrt{\xi+p+q-1}}\right)d\xi \\
&=& -\frac12\int_{\xi}^{+\infty} \frac{1}{\xi-1} - \frac{1}{p(\xi+p)} - \frac{1}{q(\xi+q)} - \frac{\sqrt{p+q}}{\sqrt{\xi + p + q - 1}}\left( \frac{1}{\xi-1} - \frac{2(p+q-1)}{(p+q)\xi} + \frac{1}{p^2(\xi+p)} + \frac{1}{q^2(\xi+q)} \right)d\xi \\
&=& \log(\sqrt{\xi + p + q - 1}+\sqrt{p+q})+\frac{\sqrt{p+q-1}}{\sqrt{p+q}}\log\xi-\frac{2\sqrt{p+q-1}}{\sqrt{p+q}}\log(\sqrt{\xi + p + q - 1}+\sqrt{p+q-1}) \\ && +\frac{1}{p}\log\frac{\sqrt{\xi + p + q - 1}+\sqrt{q-1}}{\xi+p}+\frac{1}{q}\log\frac{\sqrt{\xi + p + q - 1}+\sqrt{p-1}}{\xi+q}.
\end{eqnarray*}
It follows that
\begin{eqnarray} \label{tau_c_extreme}
t &=& (\sqrt{\xi + p + q - 1}+\sqrt{p+q})\left(\frac{\xi}{(\sqrt{\xi + p + q - 1}+\sqrt{p+q-1})^2}\right)^{\frac{\sqrt{p+q-1}}{\sqrt{p+q}}} \cdot \nonumber\\ && \cdot\left(\frac{\sqrt{\xi + p + q - 1}+\sqrt{q-1}}{\xi+p}\right)^{\frac{1}{p}}\left(\frac{\sqrt{\xi + p + q - 1}+\sqrt{p-1}}{\xi+q}\right)^{\frac{1}{q}}.
\end{eqnarray}
Relation \eqref{expphi_exp} yields
\begin{equation} \label{phi_c_extreme}
\phi = \log\left(1+\frac{\xi\sqrt{\xi + p+q-1}}{\sqrt{p+q}}\right) - \log t.
\end{equation}
Relations \eqref{tau_c_extreme} and \eqref{phi_c_extreme} define a parametrization of the solution $\phi(t)$, $t \in (0,1)$, with $\xi \in (0,+\infty)$.

By \eqref{general_immersion} the affine sphere given by the level surface $F = 0$ can be represented by the hypersurface immersion
\[ \mathbb R_{++} \times \mathbb R \ni \begin{pmatrix} \xi \\ \mu \end{pmatrix} \mapsto \begin{pmatrix} x \\ y \\ z \end{pmatrix} = \left(1+\frac{\xi\sqrt{\xi + p+q-1}}{\sqrt{p+q}}\right)^{\frac13} \begin{pmatrix} t^{-\frac13}e^{\frac{q+1}{3q}\mu} \\ t^{-\frac13}e^{-\frac{p+1}{3p}\mu} \\ t^{\frac23}e^{-\frac{\mu(p-q)}{3(p+q)}} \end{pmatrix}.
\]
Here $t$ is a function of $\xi$ defined by \eqref{tau_c_extreme}. This affine sphere is asymptotic to the cone given in case 5 of Theorem \ref{th_class}.

For the special case $p = q = 2$ we get
\begin{eqnarray*}
t &=& \left(\frac{\sqrt{\xi}}{\sqrt{\xi + 3}+\sqrt{3}}\right)^{\sqrt{3}} \cdot \frac{(\sqrt{\xi + 3}+2)(\sqrt{\xi + 3}+1)}{\xi+2}, \\
\phi &=& \log\left(1+\frac{\xi\sqrt{\xi+3}}{2}\right) - \frac{\sqrt{3}}{2}\log\xi + \sqrt{3}\log(\sqrt{\xi + 3}+\sqrt{3}) - \log\frac{(\sqrt{\xi + 3}+2)(\sqrt{\xi + 3}+1)}{\xi+2},
\end{eqnarray*}
\[ \mathbb R_{++} \times \mathbb R \ni \begin{pmatrix} \xi \\ \mu \end{pmatrix} \mapsto \begin{pmatrix} x \\ y \\ z \end{pmatrix} = \left(1+\frac{\xi\sqrt{\xi + 3}}{2}\right)^{\frac13}\begin{pmatrix} t^{-\frac13}e^{\frac{\mu}{2}} \\ t^{-\frac13}e^{-\frac{\mu}{2}} \\ t^{\frac23} \end{pmatrix}.
\]
This affine sphere is asymptotic to the boundary of the cone obtained by the homogenization of a half-disc.

$\bm{c = -2(p+q),\beta = +\infty}$. In this case $\alpha$ is positive, and we may normalize it to 1. For $t > 0$ we have the general solution
\begin{eqnarray*}
\log t &=& \int\left( \frac{1}{\xi-1} - \frac{1}{p(\xi+p)} - \frac{1}{q(\xi+q)} \right)\left(\frac12+\frac{\sqrt{p+q}}{2\xi\sqrt{\xi+p+q-1}}\right)d\xi \\ &=& \frac12\int \frac{1}{\xi-1} - \frac{1}{p(\xi+p)} - \frac{1}{q(\xi+q)}  + \frac{\sqrt{p+q}}{\sqrt{\xi + p + q - 1}}\left( \frac{1}{\xi-1} - \frac{2(p+q-1)}{(p+q)\xi} + \frac{1}{p^2(\xi+p)} + \frac{1}{q^2(\xi+q)} \right)d\xi \\ &=& \log\frac{1-\xi}{\sqrt{\xi + p + q - 1}+\sqrt{p+q}} - \frac{\sqrt{p+q-1}}{\sqrt{p+q}}\log\frac{\xi}{(\sqrt{\xi + p + q - 1}+\sqrt{p+q-1})^2} \\ && - \frac1p\log(\sqrt{\xi + p + q - 1}+\sqrt{q-1}) - \frac1q\log(\sqrt{\xi + p + q - 1}+\sqrt{p-1}) + const
\end{eqnarray*}
with $\xi \in (0,1)$. For $t < 0$ we have by virtue of \eqref{tau_positive_int}
\begin{eqnarray*}
\log|t| &=& -\int_{\xi}^{+\infty} \left( \frac{1}{\xi-1} - \frac{1}{p(\xi+p)} - \frac{1}{q(\xi+q)} \right)\left(\frac12+\frac{\sqrt{p+q}}{2\xi\sqrt{\xi+p+q-1}}\right)d\xi \\ &=& -\frac12\int_{\xi}^{+\infty} \frac{1}{\xi-1} - \frac{1}{p(\xi+p)} - \frac{1}{q(\xi+q)} + \frac{\sqrt{p+q}}{\sqrt{\xi + p + q - 1}}\left( \frac{1}{\xi-1} - \frac{2(p+q-1)}{(p+q)\xi} + \frac{1}{p^2(\xi+p)} + \frac{1}{q^2(\xi+q)} \right)d\xi \\
&=& \log\frac{\xi-1}{\sqrt{\xi + p + q - 1}+\sqrt{p+q}} - \frac{\sqrt{p+q-1}}{\sqrt{p+q}}\log\frac{\xi}{(\sqrt{\xi + p + q - 1}+\sqrt{p+q-1})^2} \\ && - \frac1p\log(\sqrt{\xi + p + q - 1}+\sqrt{q-1}) - \frac1q\log(\sqrt{\xi + p + q - 1}+\sqrt{p-1}),
\end{eqnarray*}
with $\xi \in (1,+\infty)$. The integration constant has to be chosen such that $\phi$ is smooth at $t = 0$. The two branches can be combined, yielding
\begin{eqnarray} \label{tau_c_extreme2}
t &=& \frac{1-\xi}{\sqrt{\xi + p + q - 1}+\sqrt{p+q}}\left(\frac{\sqrt{\xi + p + q - 1}+\sqrt{p+q-1}}{\sqrt{\xi + p + q - 1}-\sqrt{p+q-1}}\right)^{\frac{\sqrt{p+q-1}}{\sqrt{p+q}}}\cdot \nonumber\\ && \cdot(\sqrt{\xi + p + q - 1}+\sqrt{q-1})^{-\frac1p}(\sqrt{\xi + p + q - 1}+\sqrt{p-1})^{-\frac1q},
\end{eqnarray}
where $\xi \in (0,+\infty)$.

From \eqref{expphi_exp} we obtain
\begin{equation} \label{phi_c_extreme2}
\phi = \log\frac{\sqrt{p+q}-\xi\sqrt{\xi+p+q-1}}{\sqrt{p+q}\,t} = \log\frac{1+\xi+\frac{\xi^2}{p+q}}{1+\xi\sqrt{1+\frac{\xi-1}{p+q}}} + \log\frac{1-\xi}{t}.
\end{equation}
Relations \eqref{tau_c_extreme2} and \eqref{phi_c_extreme2} define a parametrization of the solution $\phi(t)$ on the interval $(-1,+\infty)$. By virtue of \eqref{general_immersion} the affine sphere given by the level surface $F = 0$ can then be represented by the hypersurface immersion
\[ \mathbb R_{++} \times \mathbb R \ni \begin{pmatrix} \xi \\ \mu \end{pmatrix} \mapsto \begin{pmatrix} x \\ y \\ z \end{pmatrix} = e^{\frac{\phi(\xi)}{3}} \begin{pmatrix} e^{\frac{q+1}{3q}\mu} \\ e^{-\frac{p+1}{3p}\mu} \\ e^{-\frac{\mu(p-q)}{3(p+q)}}t(\xi) \end{pmatrix}.
\]
This affine sphere is asymptotic to the cone given in case 3 of Theorem \ref{th_class}, up to a sign change in the second coordinate.

For the special case $p = q = 2$ we get
\begin{eqnarray*}
t &=& \frac{1-\xi}{(\sqrt{\xi + 3}+2)(\sqrt{\xi + 3}+1)}\left(\frac{\sqrt{\xi + 3}+\sqrt{3}}{\sqrt{\xi + 3}-\sqrt{3}}\right)^{\frac{\sqrt{3}}{2}} \\
\phi &=& \log\frac{1+\xi+\frac{\xi^2}{4}}{1+\xi\sqrt{1+\frac{\xi-1}{4}}} + \log(\sqrt{\xi + 3}+2) + \log(\sqrt{\xi + 3}+1) + \frac{\sqrt{3}}{2}\log\frac{\sqrt{\xi + 3}-\sqrt{3}}{\sqrt{\xi + 3}+\sqrt{3}},
\end{eqnarray*}
\[ \mathbb R_{++} \times \mathbb R \ni \begin{pmatrix} \xi \\ \mu \end{pmatrix} \mapsto \begin{pmatrix} x \\ y \\ z \end{pmatrix} = e^{\frac{\phi(\xi)}{3}}\begin{pmatrix} e^{\frac{\mu}{2}} \\ e^{-\frac{\mu}{2}} \\ t(\xi) \end{pmatrix}.
\]
This affine sphere is asymptotic to the boundary of the cone obtained by the homogenization of a semi-infinite strip capped by a half-disc.

\bibliography{affine_geometry,misc}
\bibliographystyle{plain}

\end{document}